\documentclass{birkjour}

\usepackage{graphicx}
\usepackage{color}
\usepackage{soul}
\numberwithin{equation}{section}

\newcommand{\R}{\mathbb{R}}

\newtheorem{Prop}{Proposition}[section]
\newtheorem{Thm}{Theorem}

\begin{document}

\title[An adaptive POD method for optimal control problems]{An adaptive POD approximation method for the control of advection-diffusion equations}
\thanks{The authors wish to acknowledge the support obtained by the following grants: ESF-OPTPDE  Network,  ITN - Marie Curie Grant n. 264735-SADCO and PRIN 2009 {"Metodi Innovativi per il Calcolo Scientifico"}.\\ The authors also wish to thank the CASPUR Consortium for its technical support}

\author[A. Alla]{A. Alla}
\address{%
Universit\`a degli studi di Roma "La Sapienza"\\
Piazzale Aldo Moro, 2
0010 Roma
Italy}
\email{alla@mat.uniroma1.it}
\author[M. Falcone]{M. Falcone}
\address{%
Universit\`a degli studi di Roma "La Sapienza"\\
Piazzale Aldo Moro, 2
0010 Roma
Italy}
\email{falcone@mat.uniroma1.it}


\date{January 22, 2012}

\begin{abstract}
We present an algorithm for the approximation of a finite horizon optimal control problem for advection-diffusion equations. The method is based on the coupling between an adaptive POD representation of the solution and a Dynamic Programming approximation scheme for the corresponding evolutive Hamilton-Jacobi equation. We discuss several features regarding  the adaptivity of the method, the role of  error estimate indicators to choose a time subdivision of the problem and the computation of the basis functions. Some test problems are presented to illustrate the method.
\end{abstract}

\maketitle

\section{Introduction}
The approximation of optimal control problems for evolutionary partial differential equations of parabolic and hyperbolic type  is a very challenging topic with a strong impact on industrial applications.
Although there is a large  number of papers dealing with several aspects of control problems from controllability to optimal control, the literature dealing with  the numerical approximation of such huge problems is rather limited. It is worth to note that  when dealing with optimal control problems for  parabolic equations  we can exploit the regularity of the solutions, regularity which is lacking for many hyperbolic equations. We also recall that the main tools is still given by the  Pontryagin maximum principle. This is mainly due to the fact that the discretization of partial differential equations already  involves a large number of variables so that the resulting finite dimensional optimization problem easily reaches the limits of what one can really compute. 
 The forward-backward system which describes Pontryagin's optimality condition  is certainly below that limit. However just solving that system one is using necessary conditions for optimality  so, in principle, there is no guarantee that these are optimal controls. By this approach for general nonlinear control problems we can obtain just open-loop control. One notable exception is the linear quadratic regulator problem for which we have a closed-loop solution given by the Riccati equation. This explains why  the most popular example for the control of evolutive partial differential equations is the control of the heat equation subject to a quadratic cost functional.\\
In recent years, new tools have been developed to deal with optimal control problems in infinite dimension. In particular, new techniques emerged to reduce the number of dimensions in the description of the dynamical system or, more in general, of the solution of the problem that one is trying to optimize. These methods are generally called {\em reduced-order methods} and include  for example the POD (Proper Orthogonal Decomposition) method  and reduced basis approximation (see \cite{PR06}). The general  idea for all this method is that, when the solution are sufficiently regular, one can represent them via Galerkin expansion so that the number of variables involved in this discretization will be strongly  reduced. In some particular case, as for the heat equation, even 5 basis functions will suffice to have a rather accurate POD representation of the solution. Having this in  mind, it is reasonable to start thinking to a different approach based on Dynamic Programming (DP) and Hamilton-Jacobi-Bellman equations (HJB). In this new approach we will first develop a reduced basis representation 
of the solution along a reference trajectory and then use this basis to set-up a control problem in the new space of coordinates. The corresponding Hamilton-Jacobi equation 
will just need 3-5 variables to represent the state of the system. Moreover, by this method one can obtain optimal control in feedback form looking at the gradient of the value function.\\
However, the solution of HJB equation it is not an easy task from the numerical point of view: the analytical solution of the HJB equation are non regular (typically, just Lipschitz continuous). 
 Optimal  control problems for ODEs were solved by Dynamic Programming, both analytically and numerically (see \cite {BCD97} for a general presentation of this theory). From the numerical point of view, this approach has been developed for many classical control problems
obtaining convergence results and a-priori error estimates (\cite{F97}, \cite{FG99} and the book \cite{FF12}).
Although this approach suffers from the curse-of-dimensionality some algorithms in high-dimension are now available (\cite{CFF04} and \cite{CCFP12}) and the coupling with POD reppresentation techniques will allow to attack by this technique optimal control problems in infinite dimension.\\
To set this paper into perspective we must say that a first tentative in this direction has been made by  Kunisch and co-authors in a series of papers \cite{KV99, KV01} for diffusion dominated equations. In particular, in the paper by Kunisch, Volkwein and Xie \cite{KVX04} one can see a feedback control approach based on coupling between POD basis approximation and HJB equations for the viscous Burgers equation. Our contribution here is twofold. The first novelty is that  we deal with advection-diffusion equations. The solutions to these equations exhibit low  regularity properties with respect to non degenerate diffusion equations so that  a rather large number of POD basis functions will be required to obtain a good approximation  if we want to compute the POD basis just once. Naturally, this increases the number of variable in the HJB approach and constitutes a is a real bottle-neck.
In order to apply the Dynamic Programming approach to this problem we have developed an adaptive technique which allows to recompute the POD basis on different sub-intervals in order to have always accurate results  without an increase of the number of basis functions. The second contribution of this paper is the way the sub-intervals are determined. In fact, we do not use a simple uniform subdivision but rather decide to recompute the POD basis when an error indicator (detailed in Section 4) is beyond a given threshold. As we will show in the sequel, this procedure  seems to be rather efficient and accurate to deal with these large scale problems. 

\section{The POD approximation method for evolutive PDEs}
We briefly describe some important features of the POD approximation, more details as well as precise results can be found in the notes 
by Volkwein  \cite{V07}. Let us consider a matrix $Y\in\R^{m\times n},$ with rank $d\leq\min\{m,n\}.$ We will call $y_j$ the $j-$th column of the matrix $Y.$ We are looking for an orthonormal basis $\{\psi_i\}_{i=1}^\ell\in\R^m$ with $\ell\leq n$ such that the minimum of the following functional is reached:
\begin{equation}\label{min_pod}
J(\psi_1,\ldots,\psi_\ell)=\sum_{j=1}^n \left\|y_j-\sum_{i=1}^\ell \langle y_j,\psi_i\rangle\psi_i \right\|^2. 
\end{equation}
The solution of this minimization problem is given in the following theorem 
\begin{Thm}\label{the_pod}
Let $Y=[y_1,\ldots,y_n]\in\R^{m\times n}$ be a given matrix with rank $d\leq\min\{m,n\}.$ Further, let $Y=\Psi\Sigma V^T$ be the Singular Value Decomposition (SVD) of $Y$, where $\Psi=[\psi_1,\ldots,\psi_m]\in\R^{m\times m}$, $V=[v_1,\ldots,v_n]\in\R^{n\times n}$ are orthogonal matrices and the matrix $\Sigma\in\R^{m\times n}$ is diagonal, $\Sigma=diag\{\sigma_1,\dots,\sigma_m\}$. Then, for any $\ell\in\{1,\ldots,d\}$ the solution to (\ref{min_pod}) is given by the left singular vectors $\{\psi_i\}_{i=1}^\ell$, i.e, by the first $\ell$ columns of $\Psi$.
\end{Thm}
\noindent
We will call the vectors $\{\psi_i\}_{i=1}^\ell$ {\em POD basis} of rank $\ell.$
\noindent
This idea is really usefull, in fact we get a solution solving an equation whose dimension is decreased with respect to the initial one. Whenever it's possible to compute a POD basis of rank $\ell,$ we get a problem with much smaller dimension of the starting one due to the fact $\ell$ is properly chosen very small.\\
\noindent
Let us consider the following ODEs system
\begin{equation}\label{PODE}
\left\{
\begin{array}{ll}
\dot{y}(s)=Ay(s)+f(s,y(s)),\;s\in(0,T]\\\\
y(0)=y_0
\end{array}\right.
\end{equation}
where $y_0\in\R^m, A\in\R^{m\times m}$ and $f:[0,T]\times\R^m\rightarrow\R^m$ is continuous and locally Lipschitz to ensure uniqueness.\\
The system (\ref{PODE}) can  be also interpreted as a semidiscrete problem, where the matrix $A$ represents the discretization in space of an elliptic operator, say Laplacian for instance. To compute the POD basis functions, first of all we have to construct a time grid $0\leq t_1\leq\ldots\leq t_n=T$ and we suppose to know the solution of (\ref{PODE}) at given time $t_j$, $j=1,\dots, N$.  We call {\em snapshots} the solution at those fixed times. For the moment we will not deal with the problem of selecting the snapshots sequence which is a difficult problem in itself, we refer the interested  readers  to \cite{KV02}). As soon as we get the snapshots sequence, by Theorem \ref{the_pod}, we will be able to compute our POD basis, namely, $\{\psi_j\}_{j=1}^\ell$.\\
Let us suppose we can write the solution in reduced form as
$$y^\ell(s)=\sum_{j=1}^\ell y_j^\ell(s)\psi_j=\sum_{j=1}^\ell\langle y^\ell(s),\psi_j\rangle\psi_j,\qquad \forall s\in[0,T]$$
substituting this formula into  (\ref{PODE}) we obtain the reduced dynamics
\begin{equation}\label{P22}
\left\{\begin{array}{ll}
\sum\limits_{j=1}^\ell\dot{y}_j^\ell (s)\psi_j=\sum\limits_{j=1}^\ell y_j^\ell(s)A\psi_j+f(s,y^\ell(s)),\qquad s\in(0,T]\\\\
\sum\limits_{j=1}^\ell y_j^\ell(0)\psi_j=y_0.
\end{array}\right.
\end{equation}
We note that our new problem (\ref{P22}) is a problem for the $\ell\leq m$ coefficient functions $y_j^\ell(s),\; j=1,\ldots,\ell.$ Thus, the problem is low dimensional and with compact notation we get:
$$\left\{\begin{array}{ll}
\dot{y}^\ell(s)=A^\ell y^\ell(s)+F(s,y^\ell(s))\\\\
y^\ell(0)=y_0^\ell\end{array}\right.$$
where
$$A^\ell\in\R^{\ell\times \ell}\qquad \mbox{with }(A^\ell)_{ij}=\langle A\psi_i,\psi_j \rangle,$$
\begin{equation*}
y^\ell=\left(\begin{array}{ccc}
y_1^\ell\\
\vdots\\
y_\ell^\ell \\
\end{array}\right):[0,T]\rightarrow\R^\ell
\end{equation*}
$F=(F_1,\ldots,F_\ell)^T:[0,T]\times\R^\ell\rightarrow\R^\ell,$
$$F_i(s,y)=\left\langle f\left(s,\sum_{j=1}^\ell y_j\psi_j\right),\psi_i\right \rangle \;\; \mbox{for } s\in[0,T]\;\; y=(y_1,\ldots y_\ell)\in\R^\ell,$$
finally obtaining the representation of $y_0$ in $\R^\ell$
\begin{equation*}
y_0^\ell=\left(\begin{array}{ccc}
\langle y_0,\psi_1 \rangle\\
\vdots\\
\langle y_0,\psi_\ell \rangle \\
\end{array}\right)\in\R^\ell.
\end{equation*}
In order to apply the POD method to our optimal control problem, the number  $\ell$ of POD basis functions is crucial. In particular we would like to keep $\ell$ as low as possible still capturing the behaviour of the original dynamics. The problem is to define an indicator of the accuracy of our POD approximation. A good choice for this indicator  is the following ratio
\begin{equation}\label{ind:ratio}
\mathcal{E}(\ell)=\dfrac{\sum\limits_{i=1}^\ell \sigma_i}{\sum\limits_{i=1}^d \sigma_i}.
\end{equation}
where the $\sigma_i$ are the the singular value obtained by the SVD.

As much $\mathcal{E}(\ell)$ is close to one as much our approximation will be improved. This is strictly related to the truncation error due to the projection of $y_j$ onto the space generated by the orthonormal basis $\{\psi\}_{i=1}^\ell,$ in fact:
\[
J(\psi_1,\ldots,\psi_\ell)=\sum_{j=1}^n \left\|y_j-\sum_{i=1}^\ell \langle y_j,\psi_i \rangle\psi_i \right\|^2=\sum_{i=\ell+1}^d \sigma_i^2
\]

\section{An optimal control problem}
We will present this approach for the finite horizon control problem. Consider the controlled system
\begin{equation}\label{eq:consys}
\left\{ \begin{array}{l}
\dot{y}(s)=f(y(s),u(s),s), \;\; s\in(t,T]\\
y(t)=x\in\R^n,
\end{array} \right.
\end{equation}
we will denote by $y:[t,T]\rightarrow\R^n$ its the solution, by $u$ the  control $u:[t,T]\rightarrow\R^m$, $f:\R^n\times\R^m\rightarrow\R^n$, $s\in(t,T]$ and by
\[\mathcal{U}=\{u:[0,T]\rightarrow U \}
\]
the set of admissible controls where $U\subset \R^m$ is a compact set.  Whenever we want to emphasize the depence of the solution from the control $u$ we will write $y(t;u)$. Assume that there exists a unique solution trajectory for \eqref{eq:consys} provided the controls are measurable (a precise statement can be found in \cite{BCD97}). For  the finite horizon  optimal control problem the cost functional will be given by 
\begin{equation}
\min_{u\in\mathcal{U}} J_{x,t}(u):=\int_t^T L(y(s,u),u(s),s)e^{-\lambda s}\, ds+g(y(T))
\end{equation}
where $L:\R^n\times\R^m\rightarrow\R$ is the running cost and $\lambda\geq0$ is the discount factor.\\
The goal is to find a state-feedback control law $u(t)=\Phi(y(t),t),$ in terms of the state equation $y(t),$ where $\Phi$ is the feedback map. To derive optimality conditions we use the well-known {\em dynamic programming principle} due to Bellman (see \cite{BCD97}). We first define the value function:
\begin{equation}
v(x,t):=\inf\limits_{u\in\mathcal{U}} J_{x,t}(u)
\end{equation}
\begin{Prop}[DPP]
For all $x\in\R^n$and $0\leq\tau\leq t$ then:
\begin{equation}\label{dpp}
v(x,t)=\min_{u\in\mathcal{U}}\left\{\int_t^\tau L(y(s),u(s),s) e^{-\lambda s}\;ds+ v(y,t-\tau)\right\}.
\end{equation}
\end{Prop}
\noindent
Due to (\ref{dpp}) we can derive the {\em Hamilton-Jacobi-Bellman} equations (HJB):
\begin{equation}\label{HJB}
-\dfrac{\partial v}{\partial t}(y,t)=\min_{u\in U }\left\{L(y,u,t)+\nabla v(y,t) \cdot f(y,u,t)\right\}.
\end{equation}
This is nonlinear partial differential equation of the first order which is hard to solve analitically although a general theory of weak solutions is available  \cite{BCD97}.  Rather we can solve it numerically by means of a finite differences or semi-Lagrangian schemes (see the book \cite{FF12} for a comprehensive analysis of approximation schemes for Hamilton-Jacobi equations). For a semi-Lagrangian discretization one starts by a discrete version of (HJB) by discretizing the underlined control problem and then project the semi-discrete scheme on a grid obtaining the  fully discrete scheme
$$\left\{\begin{array}{ll}
v_i^{n+1}=\min\limits_{u\in U}[\Delta t\,L(x_i,n\Delta t,u)+I[v^n](x_i+\Delta t\, F(x_i,t_n,u))]\\\\
v_i^0=g(x_i).
\end{array}\right.$$
with $x_i=i\Delta x,\;t_n=n\Delta t,\;v^n_i:=v(x_i,t_n)$ and $I[\cdot]$ is an interpolation operator which is necessary to compute the value of $v^n$ at the point $x_i+\Delta t\, F(x_i,t_n,u)$ 
(in general, this point will not be a node of the grid). The interested reader will find in \cite{FG99} a detailed presentation of the scheme and a priori error estimates for its numerical approximation. 

Note that, we also need to compute the minimum in order to get the value $v_i^{n+1}$. Since $v^n$ is not a smooth function, we compute the minimum by means of a minimization method which does not use derivatives (this can be done by the Brent algorithm as in \cite{CFF04}).



As we already told the HJB allows to compute the optimal feedback via the value function, but there are two major difficulties: the solution of an HJB equation are in general non-smooth and the approximation in high dimension is not feasible. The request to solve an HJB in high dimension comes up naturally whenever we want to control evolutive PDEs. Just to give an idea, if we build a grid in $[0,1]\times[0,1]$ with a discrete step $\Delta x=0.01$ we have $10^4$ nodes: to solve an HJB in that dimension is simply impossible. Fortunatelly, the POD method allows us to obtain reduced models even for complex dynamics. 
Let us focus on the following abstract problem: 
\begin{equation}\label{pabs}
\left\{\begin{array}{ll}
\dfrac{d}{ds} \langle y(s),\varphi \rangle_H+a(y(s),\varphi)=\langle B(u(s),\varphi\rangle_{V',V}\quad \forall \varphi\in V\\\\
y(t)=y_0\;\;\in H,
\end{array}\right.
\end{equation}
where $B:U\rightarrow V'$ is a linear and continuous operator. We assume that a space of admissible controls $\mathcal{U}_{ad}$ is given in such a way that  for each $u\in\mathcal{U}_{ad}$ and $y_0\in H$ there exists a unique solution $y$ of (\ref{pabs}).
$V$ and $H$ are two Hilbert spaces, with $\langle\cdot,\cdot\rangle_H$ we denote the scalar product in  $H;$ $a:V\times V\rightarrow \R:$ is symmetric coercive and bilinear.
Then, we introduce the cost functional of the finite horizon problem
$$\mathcal{J}_{y_0,t}(u):=\int_t^T L(y(s),u(s),s) e^{-\lambda s} \;ds+g(y(T)),$$
where $L:V\times U\times[0,T]\rightarrow \R.$ The optimal control problem is 
\begin{eqnarray}\label{KP}
&&\min\limits_{u\in\mathcal{U}_{ad}} \mathcal{J}_{y_0,t}(u)\\
\hbox{subject to the constraint: }&& y\in W_{loc}(0,T;V)\times\mathcal{U}\hbox{ solves }(\ref{pabs})\nonumber
\end{eqnarray}
with $W_{loc}(0,T)=\bigcap_{T>0}W(0,T),$ where $W(0,T)$ is the standard Sobolev space:
$$W(0,T)=\{\varphi\in L^2(0,T;V), \varphi_t\in L^2(0,T;V')\}.$$
\medskip

\noindent
The model reduction approach for an optimal control problem (\ref{KP}) is based on the Galerkin approximation of dynamic with some informations on the controlled dynamic (snapshots).
To compute a POD solution for (\ref{KP}) we make the following ansatz
\begin{equation}\label{K314}
y^\ell(x,s)=\sum_{i=1}^\ell w_i(s)\psi_i(x).
\end{equation}
where $\{\psi\}_{i=1}^\ell$ is the POD basis computed as in the previous section.\\
\noindent
We introduce mass and stiffness matrix:
$$M=((m_{ij}))\in\R^{\ell\times \ell}\mbox{ with } m_{ij}=\langle\psi_j,\psi_i\rangle_H,$$
$$S=((s_{ij}))\in\R^{\ell\times \ell}\mbox{ with } m_{ij}=a(\psi_j,\psi_i),$$
and the control map $b:U\rightarrow \R^\ell$ is defined by:
$$u\rightarrow b(u)=(b(u)_i)\in\R^\ell\mbox{ with } b(u)_i=\langle Bu,\psi_i\rangle_H.$$
The coefficients of the initial condition $y^\ell(0)\in\R^\ell$ are determined by $w_i(0)=(w_0)_i=\langle y_0,\psi\rangle_X, \;\;1\leq i\leq \ell,$ and the solution of the reduced dynamic problem is denoted by $w^\ell(s)\in\R^\ell.$ Then,  the Galerkin approximation is given by
\begin{equation}\label{KPl}
\min J^\ell_{w^\ell_0,t} (u)\\
\end{equation}
with $u\in\mathcal{U}_{ad}$ and $w$ solves the following equation:
\begin{equation}\label{KPl1}
 \left\{\begin{array}{ll}
\dot {w}^\ell(s)=F(w^\ell(s),u(s),s)\;\;\;s>0,\\\\
w^\ell(0)=w_0^\ell.\end{array}\right.
\end{equation}
The cost functional is defined:
$$J^\ell_{w_0^\ell, t}(u)=\int_0^T L(w^\ell(s),u(s),s)e^{-\lambda s}\; dt+g(w^\ell(T)),$$
with  $w^{\ell}$ and $y^\ell$ linked to (\ref{K314}) and the nonlinear map $F:\R^\ell\times U\rightarrow \R^\ell$ is given by
$$F(w^\ell,u,s)=M^{-1}(-Sw^\ell(s)+b(u(s))).$$
The value function $v^\ell$, defined for the initial state $w_0\in\R^\ell,$ 
$$v^\ell(w^\ell_0,t)=\inf_{u\in\mathcal{U}_{ad}} J^\ell_{w_0^\ell, t}(u)$$ 
and $w^\ell$ solves (\ref{KPl}) with the control $u$ and initial condition $w_0.$

\medskip
\noindent
 We give an idea how we have computed the intervals for reduced HJB. HJBs are defined in $\R^n,$ but we have restricted our numerically domain $\Upsilon_h$ which is a bounded subset of $\R^n.$ This is justified since $y+\Delta t F(y,u)\in\Upsilon_h$ for each 
$y\in\Upsilon_h$ and $u\in\mathcal{U}_{ad}.$ We can chose $\Upsilon_h=[a_1,b_1]\times[a_2,b_2]\times\ldots[a_\ell,b_\ell]$ with $a_1\geq a_2\geq\ldots \geq a_\ell.$ How should we compute these intervals $[a_i,b_i]$?\\
\noindent
Ideally the intervals should be chosen so that the dynamics contains all the components of the controlled trajectory.  Moreover, they should be encapsulated because we expect that their importance should decrease monotonically  with their index  and that our interval lengths decrease quickly.\\
Let us suppose to discretize the space control $U=\{u_1,\ldots,u_M\}$ where $U$ is symmetric, to be more precise if $\bar{u}\in U\Rightarrow -\bar{u}\in U.$\\
\noindent
Hence, if $y^\ell(s)=\sum\limits_{i=1}^\ell \langle y(s),\psi_i\rangle\psi_i=\sum_{i=1}^\ell w_i(s)\psi_i,$ as a consequence, the coef\-ficients $w_i(s)\in[a_i,b_i].$ We consider the trajectories solution $y(s,u_j)$ such that the control is constant $u(s)\equiv u_j$ for each $t_j$, $j=1,\ldots, M.$ Then,  we have
$$y^\ell(s,u_j)=\sum_{i=1}^\ell \langle y(s,u_j),\psi_i\rangle\psi_i.$$
We write $y^\ell(s,u_j)$ to stress the dependence on the  constant control $u_j.$ Each trajectory $y^\ell(s,u_j)$ has some coefficients $w_i^{(j)}(t)$ for $i=1,\ldots, \ell,\, j=1,\ldots, M.$
 The coefficients $w_i^{(j)}(s)$ will belong to intervals of the type $[\underbar{w}_i^{(j)},\overline{\rm{w}}_i^{(j)}]$ where we chose for $i=1,\ldots,\ell,$ $a_i, b_i$ such that:
$$a_i\equiv\min\{\underbar{w}_i^{(1)},\ldots,\underbar{w}_i^{(M)}\}$$
$$b_i\equiv\max\{\overline{\rm{w}}_i^{(1)},\ldots,\overline{\rm{w}}_i^{(M)}\}.$$
Then, we have a method to compute the intervals and we turn our attention to the  numerical solution of an optimal control problem for evolutive equation, as we will see in the following section.


\section{Adapting POD approximation}
\begin{figure}

\hspace{-0.3cm}
\begin{minipage}{\textwidth}
\includegraphics[scale=0.3]{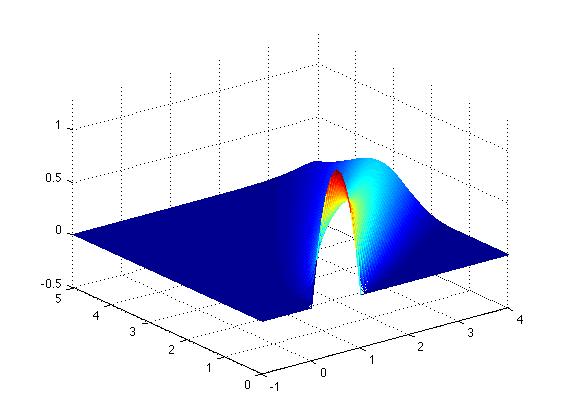}\;
\includegraphics[scale=0.3]{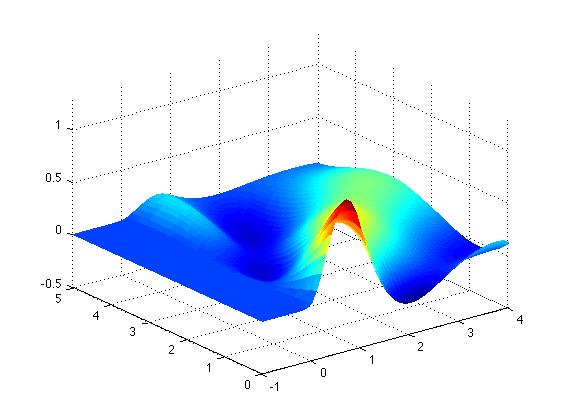}
\end{minipage}

\vspace{0.2cm}
\hspace{0.2cm}(a)\hspace{5.5cm}(b)
\hspace{-2cm}
\begin{minipage}{\textwidth}
\includegraphics[scale=0.3]{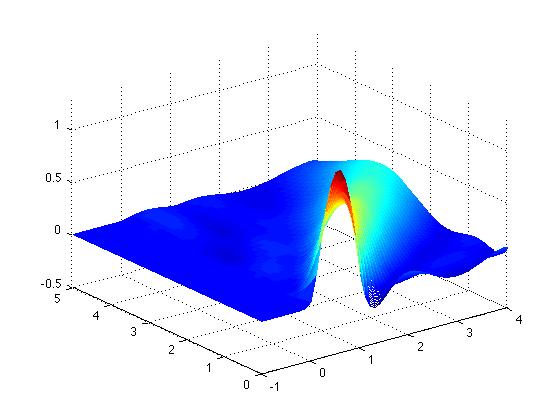}\;
\includegraphics[scale=0.3]{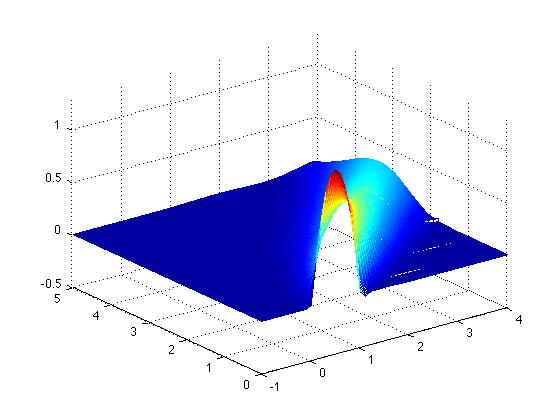}
\end{minipage}

\vspace{0.2cm}
\hspace{0.2cm}(c)\hspace{5.5cm}(d)
\caption{Equation (\ref{eq:test}):(a) solved with finite difference; (b) POD-Galerkin approximation with 3 POD basi; (c) solved via POD-Galerkin approximation with 5 POD basis; (d) Adapting 3 POD basis functions.}
\label{test10}
\end{figure}

We now present an adaptive method  to compute POD basis. Since our final goal is to obtain  the optimal feedback law by means of HJB equations, we will have a big constraint on the number of variables  in the state space for numerical solution of an HJB.\\
\noindent
We will see that, for  a parabolic equation, one can  try to solve the problem with only three/four POD basis functions; they are enough to describe the solution in a rather accurate way. In fact the singular values decay pretty soon and it's easier to work with a really low-rank dimensional problem.\\
\noindent
On the contrary, hyperbolic equations do not have this  nice property for their singular values and they will require a rather large set of POD basis functions to get accurate results. Note that we can not  follow the  approach suggested  in \cite{RV10} because we can not  add more basis functions when it turns to be necessary due to the constraint already mentioned. Then, it is quite natural to split  the problem into subproblems having different POD basis functions. The crucial point is to decide the splitting in order to have the same number of basis functions in each subdomain with a guaranteed accuracy in the approximation.\\
Let us first give an illustrative  example  for the parabolic case, considering a 1D advection-diffusion equation:
\begin{equation}\label{eq:test}
\left\{\begin{array}{ll}
y_s(x,s)-\varepsilon y_{xx}(x,s)+cy_x(x,s)=0\\
y(x,0)=y_0(x),
\end{array}\right.
\end{equation}
with $x\in[a,b], s\in[0,T], \varepsilon,c\in\R.\\$
\noindent
We use a finite difference approximation for this equation based on an explicit
 Euler method in time combined with the standard centered approximation of the second order term and   with an up-wind correction  for the advection term.  The snapshots will be taken from the sequence generated by the finite difference method. The final time is $T=5$,  moreover $a=-1$, $b=4$. The initial condition is $y_0(x)=5x-5x^2,$ when $0\leq x \leq 1$, 0 otherwise.\\
For  $\varepsilon=0.05$ and $c=1$ with only 3 POD basis functions, the approximation fails (see Figure \ref{test10}). Note that in this case the advection is dominating the diffusion, a low number of POD basis functions will not suffice to get an accurate approximation (Figure 1.b). However, the adaptive method which only uses 3 POD basis functions will give accurate results (Figure 1.d).  \\
\noindent

\noindent
The idea which is behind the adaptive method is the following: we do not consider  all the snapshots together in the whole interval $[0,T]$ but we group them. Instead of taking into account the whole interval $[0,T],$ we prefer to split it in sub-intervals 
\[ 
[0,T]=\cup_{k=0}^{K}[T_k, T_{k+1}] 
\]
 where $K$ is a-priori unknown, $T_0=0, T_{K}=T$ and $T_k=t_i$ for some $i.$ 
 In this way, choosing properly the length of the $k-$th interval $[T_k, T_{k+1}],$ we consider only the snapshots falling in that sub-interval, typically there will be at least three snapshots in every sub-interval. Then we have enough informations in every sub-interval and we can apply  the standard routines (explained in Section 2)  to get a "local"  POD basis.  \\
\noindent
Now let us explain how to divide our time interval $[0,T]$. We will choose a parameter to check the accuracy of the POD approximation and define a threshold. Above that threshold we loose in accuracy and we need to compute a new POD basis.  A good parameter to check the accuracy is $\mathcal{E}(\ell)$ (see (\ref{ind:ratio})), as it was suggested by several authors. The method to define the splitting of $[0,T]$  and the size of every sub-interval works as follows. We start computing the SVD of the matrix $Y$ that gives us informations about our dynamics in the whole time interval. We check the accuracy at every $t_i$, $i=1, \dots N$,  and  if  at $t_k$ the indicator is above the  tolerance
we set $T_1=t_k$ and  we divide the interval in two parts, $[0,T_1)$ and $(T_1, T]$. Now we just consider the snapshots related the solution up to the time $T_1$.  Then we iterate this idea until the indicator is below the threshold. When the first interval is found,  we restart the procedure in the interval $[T_1, T]$ and we stop when we reach the final time $T$. Note that the extrema of every interval coincide by construction  with one of our discrete times $t_i=i\Delta t$ so  that the global solution is easily obtained linking all the sub-problems which always have a snapshot as initial condition. A low value for the threshold will also guarantee that we will not have big jumps passing from one sub-interval to the next.


\noindent
This idea can be applied also when we have a  controlled dynamic (see (\ref{eq:test_cont})). First of all we have to decide how to collect the snapshots, since the control $u(t)$ is completely unknown. One can make a guess and use the dynamics and the functional corresponding to that guess, by these informations  we can compute the POD basis. Once the POD basis is obtained we will get the optimal feedback law after having solved a reduced HJB equation as we already explained.
Let us summarize the method in the following  step-by-step presentation.\\
\noindent
\medskip

\noindent
{\tt {\bf ALGORITHM}\\
{\em Start: } Inizialization\\
{\bf Step 1:}  collect the  snapshots in [0,T]\\
{\bf Step 2:} divide $[0,T]$ according to $\mathcal{E}(\ell)$\\
{\em For} i=0 to N-1\\
{\em Do}\\
\phantom{x}\hspace{0.2cm}{\bf Step 3:} apply SVD to get the POD basis in each sub-interval $[t_i,t_{i+1}]$\\
\phantom{x}\hspace{0.2cm}{\bf Step 4:} discretize the space of controls\\
\phantom{x}\hspace{0.2cm}{\bf Step 5:} project the dynamics onto  the (reduced) POD space\\
\phantom{x}\hspace{0.2cm}{\bf Step 6:} select the intervals for the POD reduced variables\\
\phantom{x}\hspace{0.2cm}{\bf Step 7:} solve the corresponding HJB in the reduced space\\
\phantom{x}\hspace{1cm}for the interval $[t_i,t_{i+1}]$\\
\phantom{x}\hspace{0.2cm}{\bf Step 8:} go back to the original coordinate space\\
{\em End }
}

\section{Numerical experiments}\label{tests}
In this section we present some numerical tests for the  controlled heat equation and for the advection-diffusion equation with a quadratic cost  functional.
Consider the following advection-diffusion equation:
\begin{equation}\label{eq:test_cont}
\left\{\begin{array}{ll}
y_s(x,s)-\varepsilon y_{xx}(x,s)+cy_x(x,s)=u(s)\\
y(x,0)=y_0(x),
\end{array}\right.
\end{equation}
with $x\in[a,b]$,  $s\in[0,T]$, $\varepsilon\in\R_+$ and $c\in\R.\\$Note that changing the parameters $c$ and $\varepsilon$ we can obtain the heat equation ($c=0$) and the advection equation ($\varepsilon=0$).
The functional to be minimized is
\begin{equation}\label{lqr_oss}
J_{y_0,t}(u(\cdot))=\int_0^T  ||y(x,s)-\widehat y(x,s)||^2 + R||u(s)||^2 \;ds,
\end{equation}
\noindent
i.e. we want to stay close to a reference trajectory $\widehat y$ while minimizing the norm of $u$.
Note that we dropped the discount factor setting $\lambda=0$. Typically in our test problems $\widehat y$ is obtained by applying a particular control $\widehat u$  to the dynamics. 
The numerical simulations reported in this papers  have been made on a server SUPERMICRO 8045C-3RB  with 2 cpu Intel Xeon Quad-Core 2.4 Ghz and 32 GB RAM
under SLURM ({\tt https://computing.llnl.gov/linux/slurm/}).

\begin{figure}

\hspace{-0.3cm}
\begin{minipage}{\textwidth}
\includegraphics[scale=0.4]{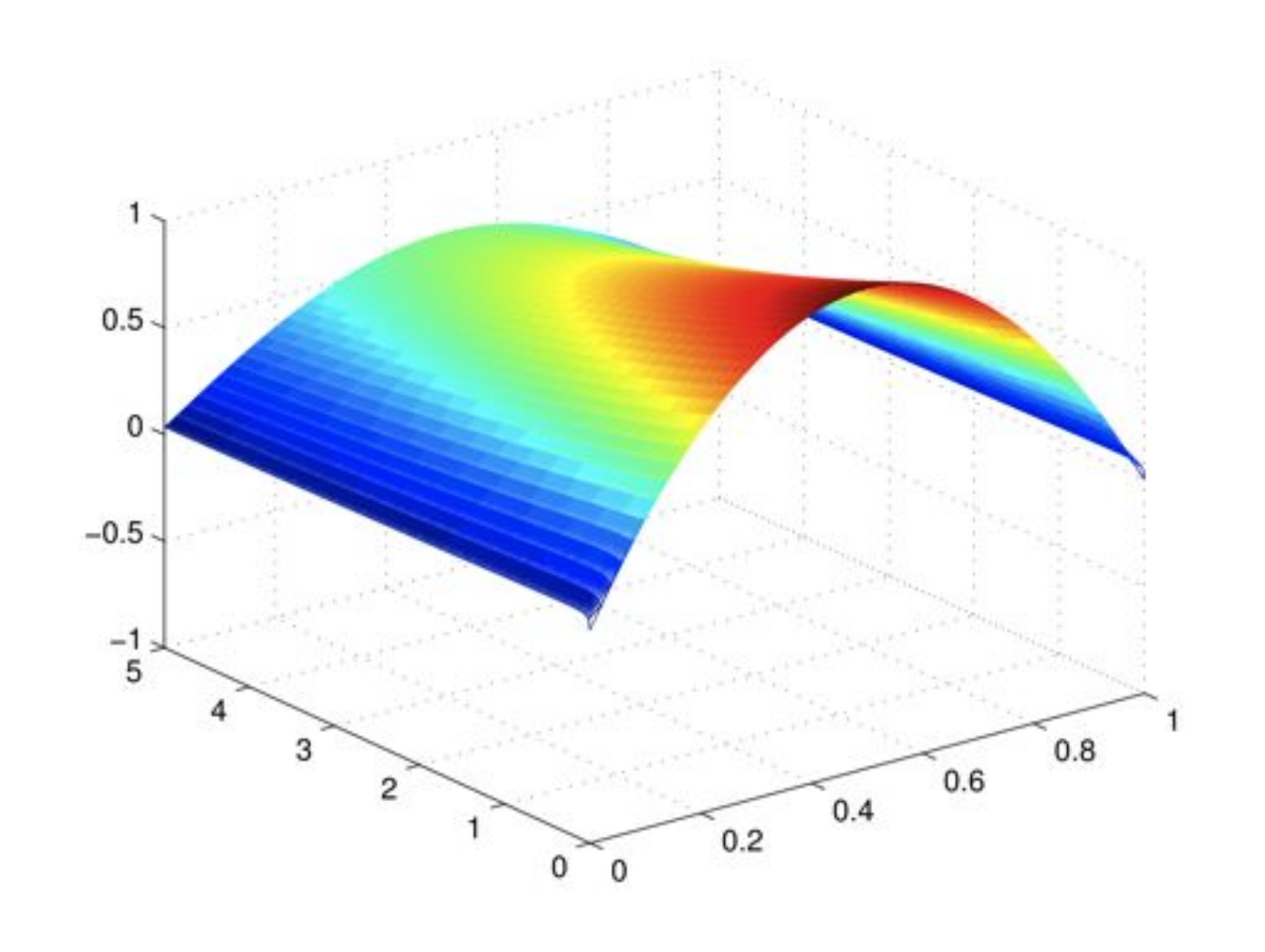}\;
\includegraphics[scale=0.4]{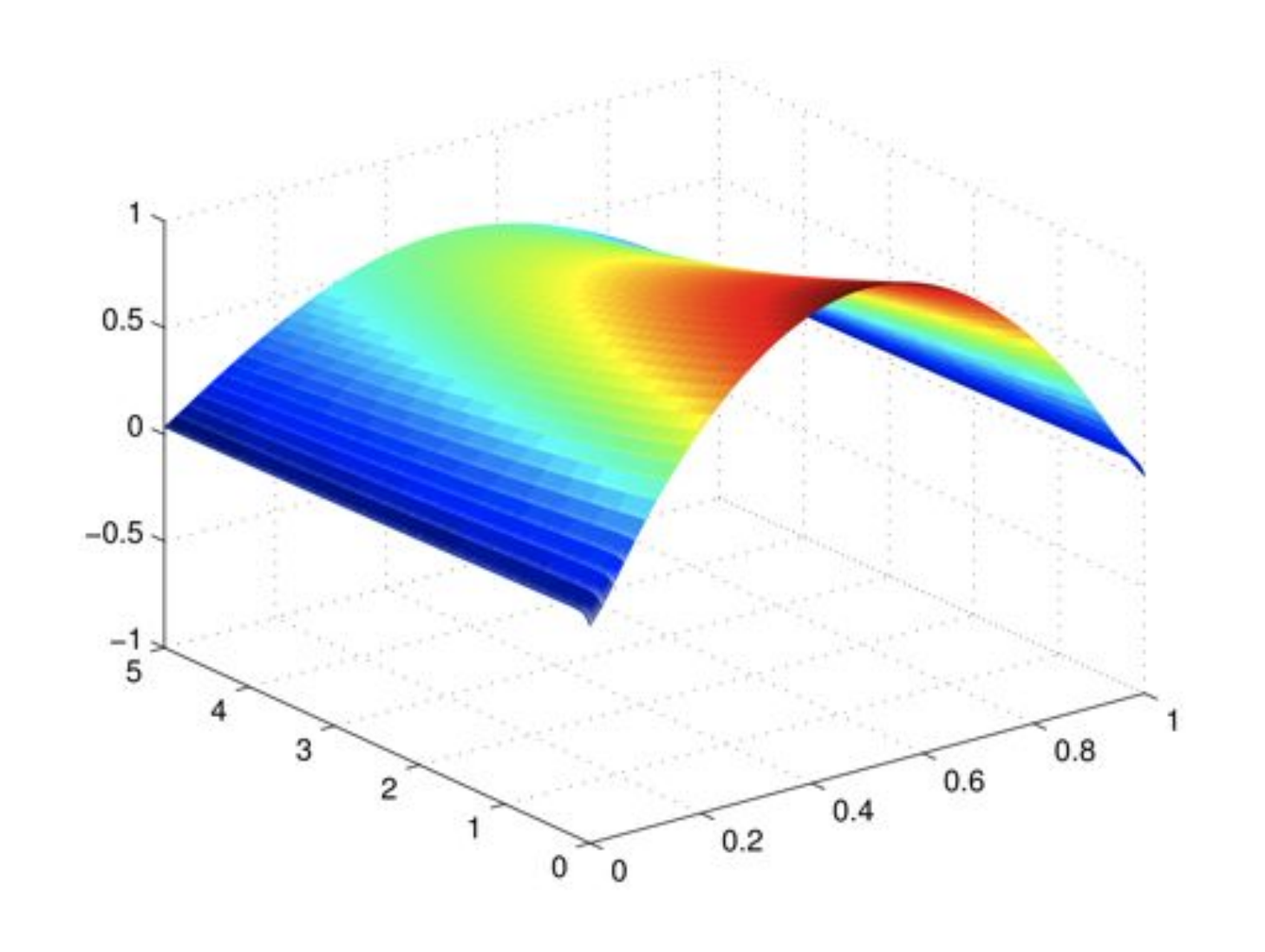}
\end{minipage}

\vspace{0.2cm}
\hspace{0.2cm}(a)\hspace{5.5cm}(b)
\hspace{-2cm}
\begin{minipage}{\textwidth}
\includegraphics[scale=0.4]{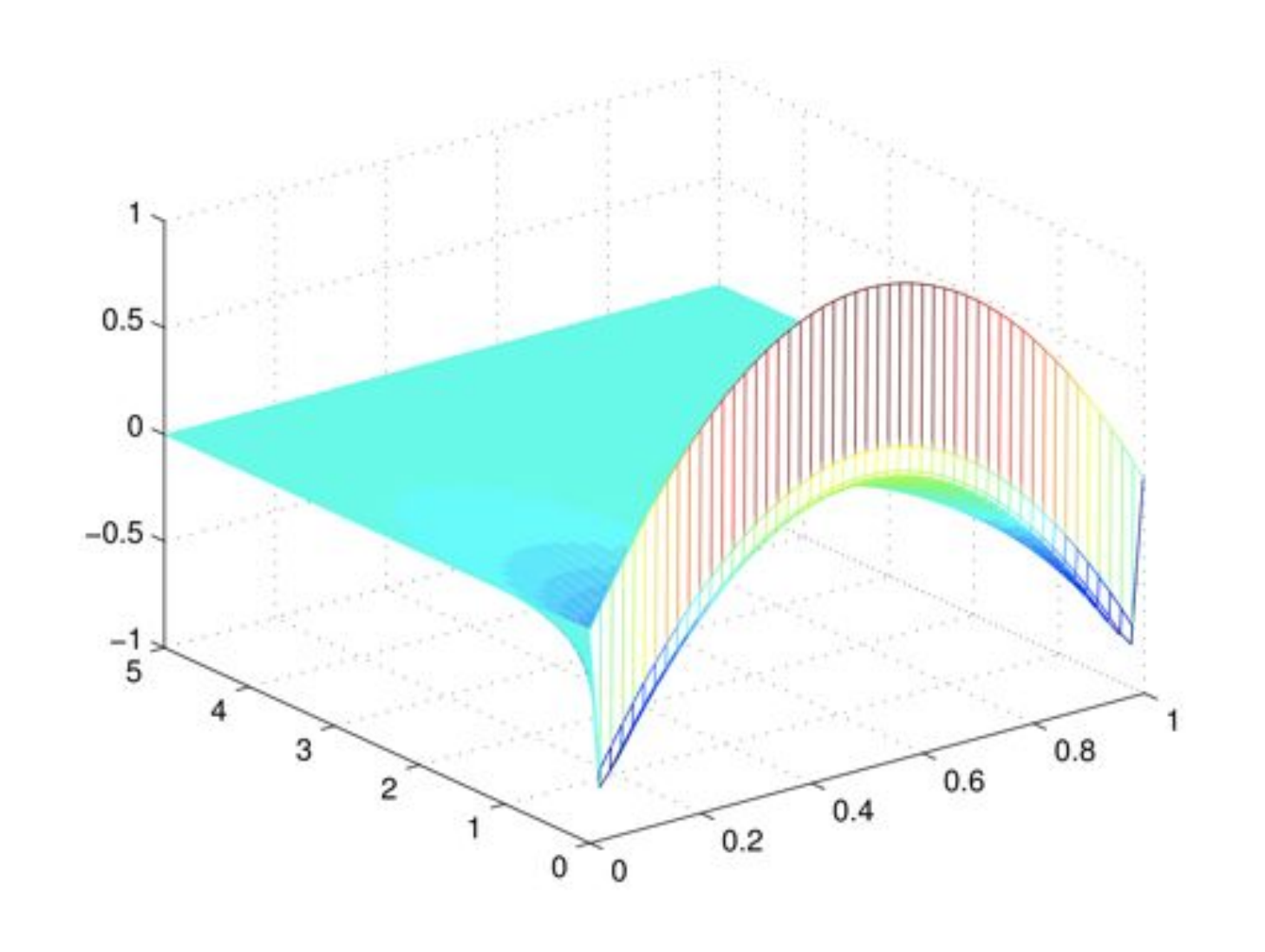}\;
\includegraphics[scale=0.4]{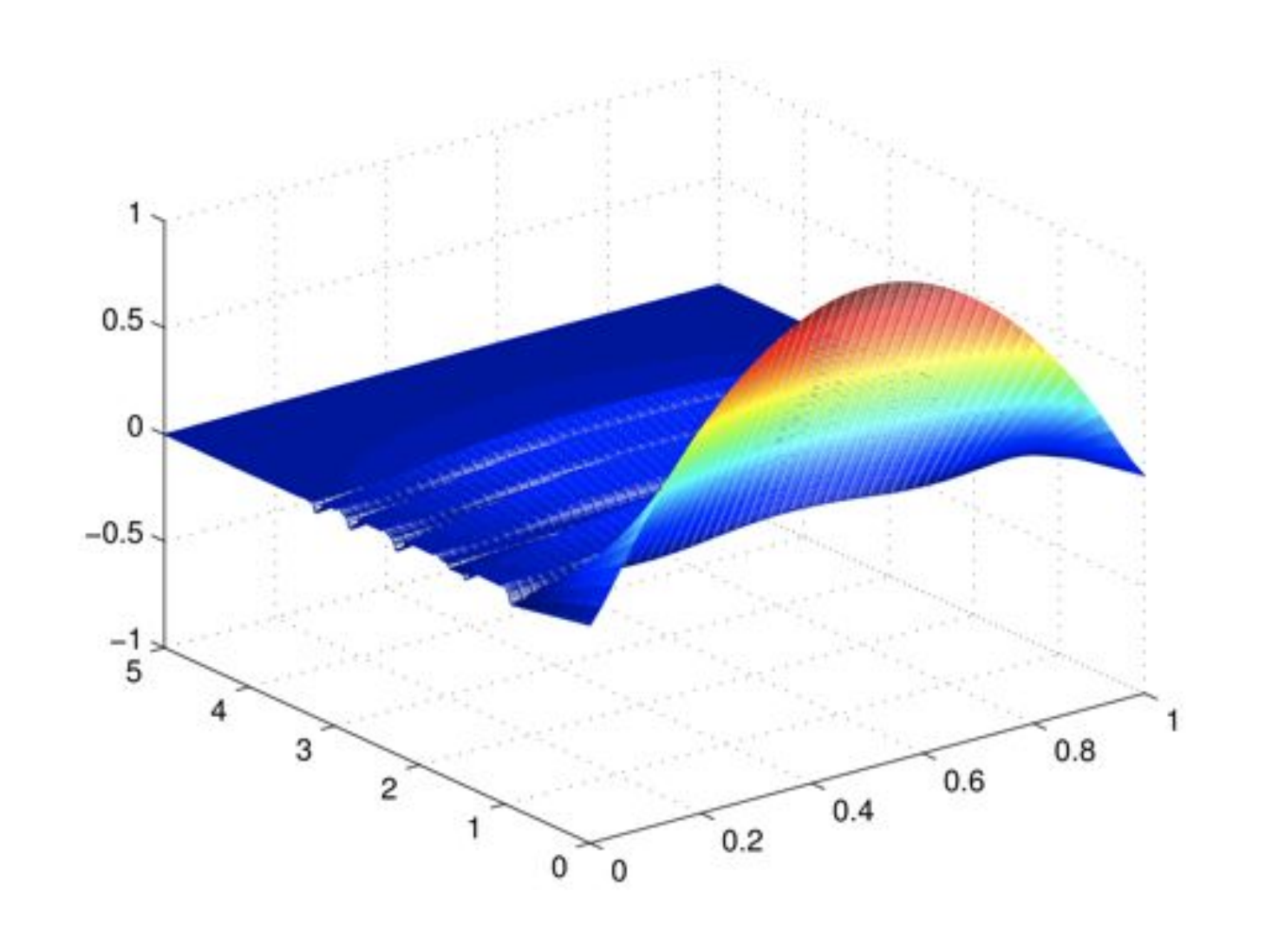}
\end{minipage}

\vspace{0.2cm}
\hspace{0.2cm}(c)\hspace{5.5cm}(d)
\caption{Test 1:(a) Heat Equation without control; (b) Heat Equation without control, 3 POD basis; (c) Controlled solution with LQR-MATLAB; (d) Approximate solution POD (3 basis functions) + HJB.}
\label{test}
\end{figure}

\medskip
\paragraph{Test 1: Heat equation with smooth initial data}
We compute the snapshots with a centered/forward Euler scheme with space step $\Delta x=0.02$, and time step $\Delta t=0.012$, $\varepsilon=1/60, c=0, R=0.01$ and $T=5$.
The initial condition is $y_0(x)=5x-5x^2,$ and $\widehat y(x,s)=0.$
In Figure \ref{test} we compare four different approximations concerning the heat equation: (a)  is the solution for $\widehat u(t)=0$, (b) is its approximation via POD (non adaptive), (c) is the direct LQR  solution computed by MATLAB without POD and, finally, the approximate optimal solution obtained coupling POD and HJB. The approximate value function is  computed for $\Delta t=0.1$ $\Delta x=0.1$  whereas the optimal trajectory as been obtained  with $\Delta t=0.01.$ Test 1, and even Test 2, have been solved in about half an hour of CPU time.\\
Note that in this  example the approximate solution is rather accurate because the regularity of the solution is high due to the diffusion term. Since in the limit  the solution tends to the average value the choice of the snapshots will not affect too much the solution, i.e.  even with a  rough choice of the snapshots will give us a good approximation.  The difference between Figure 2c and Figure 2d is due to the fact that the control space is continuous for 2c  and discrete for 2d.

\begin{figure}

\hspace{-0.3cm}
\begin{minipage}{\textwidth}
\includegraphics[scale=0.4]{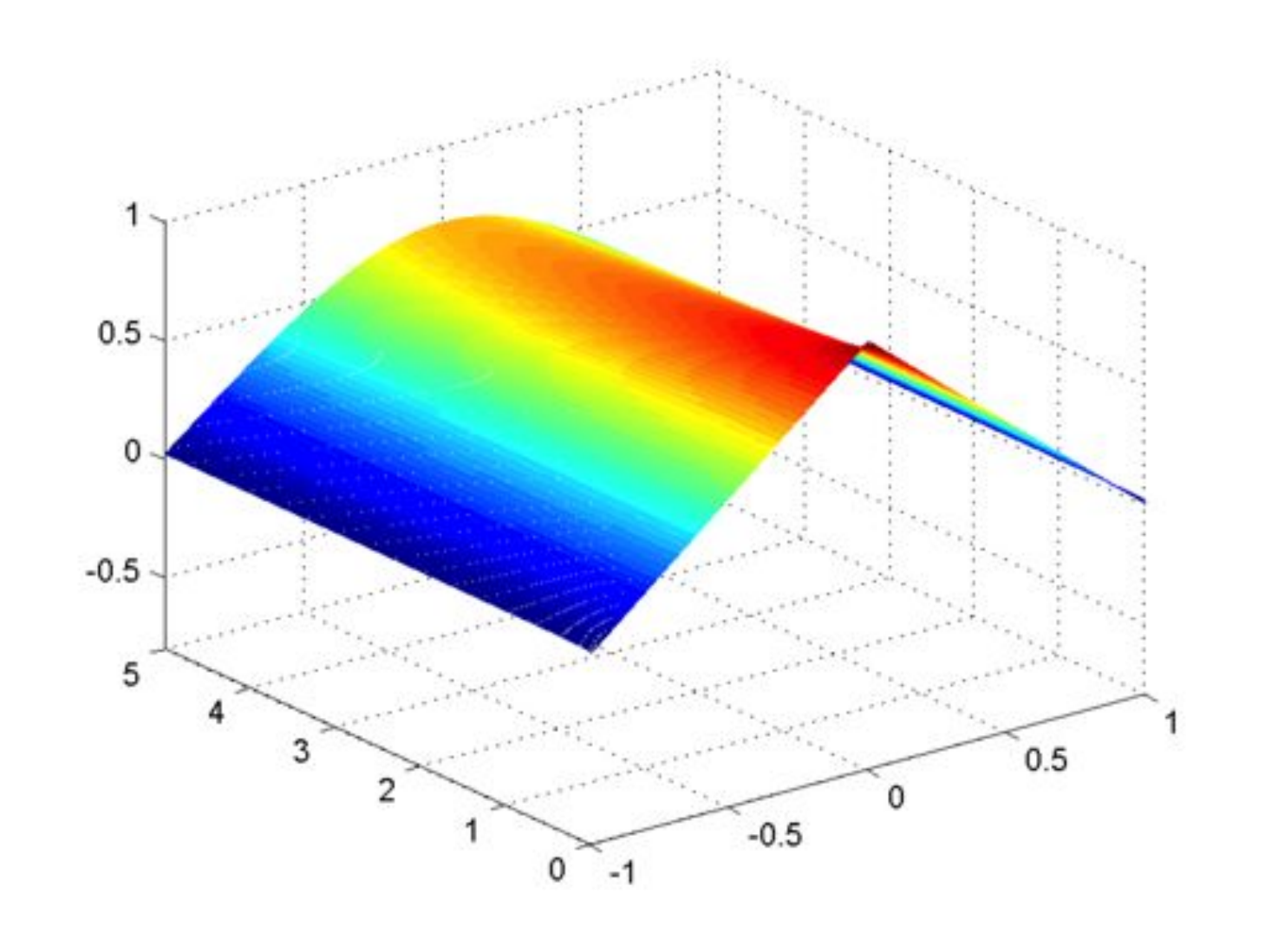}\;
\includegraphics[scale=0.4]{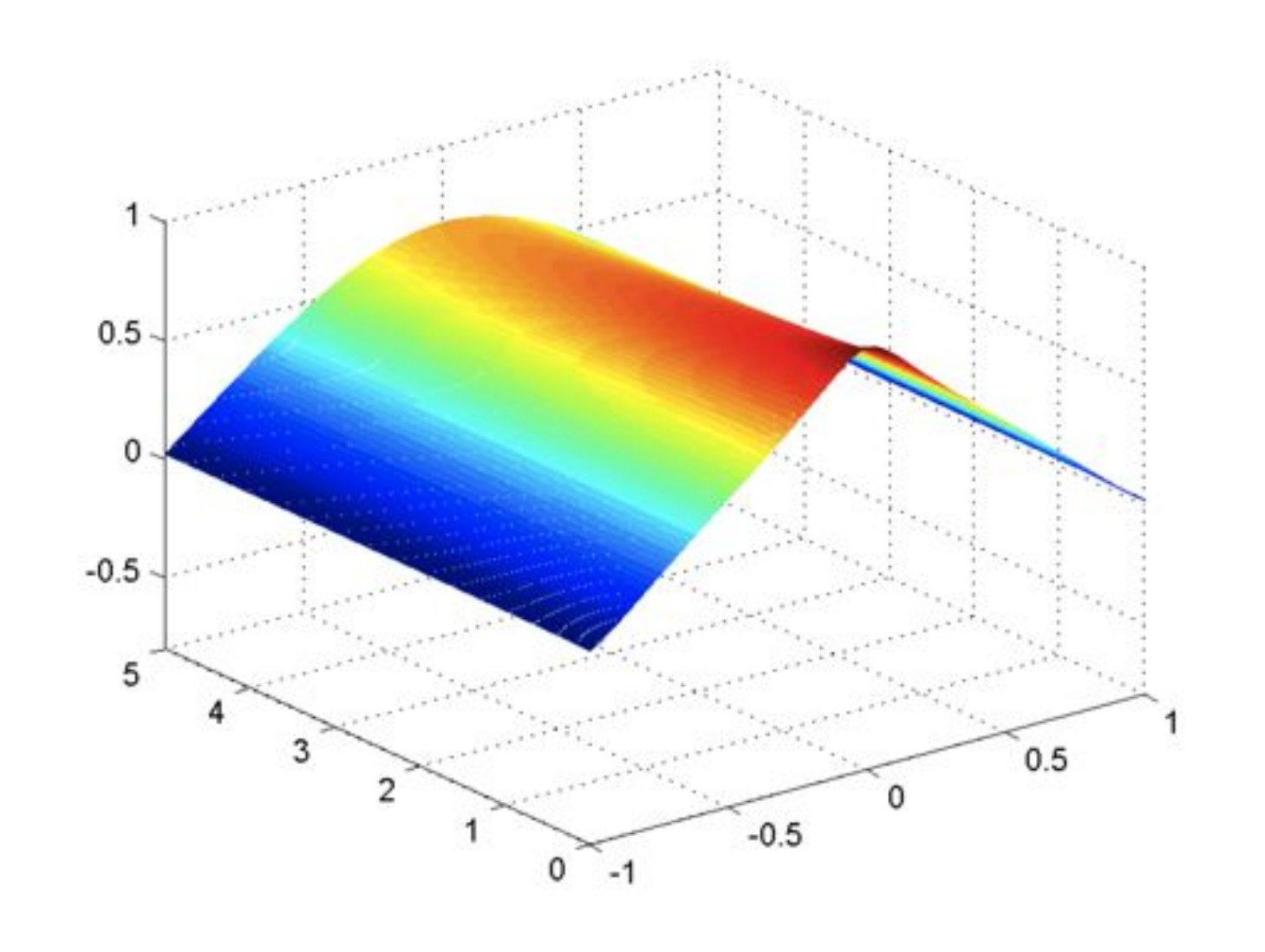}
\end{minipage}

\vspace{0.2cm}
\hspace{0.2cm}(a)\hspace{5.5cm}(b)
\hspace{-2cm}
\begin{minipage}{\textwidth}
\includegraphics[scale=0.4]{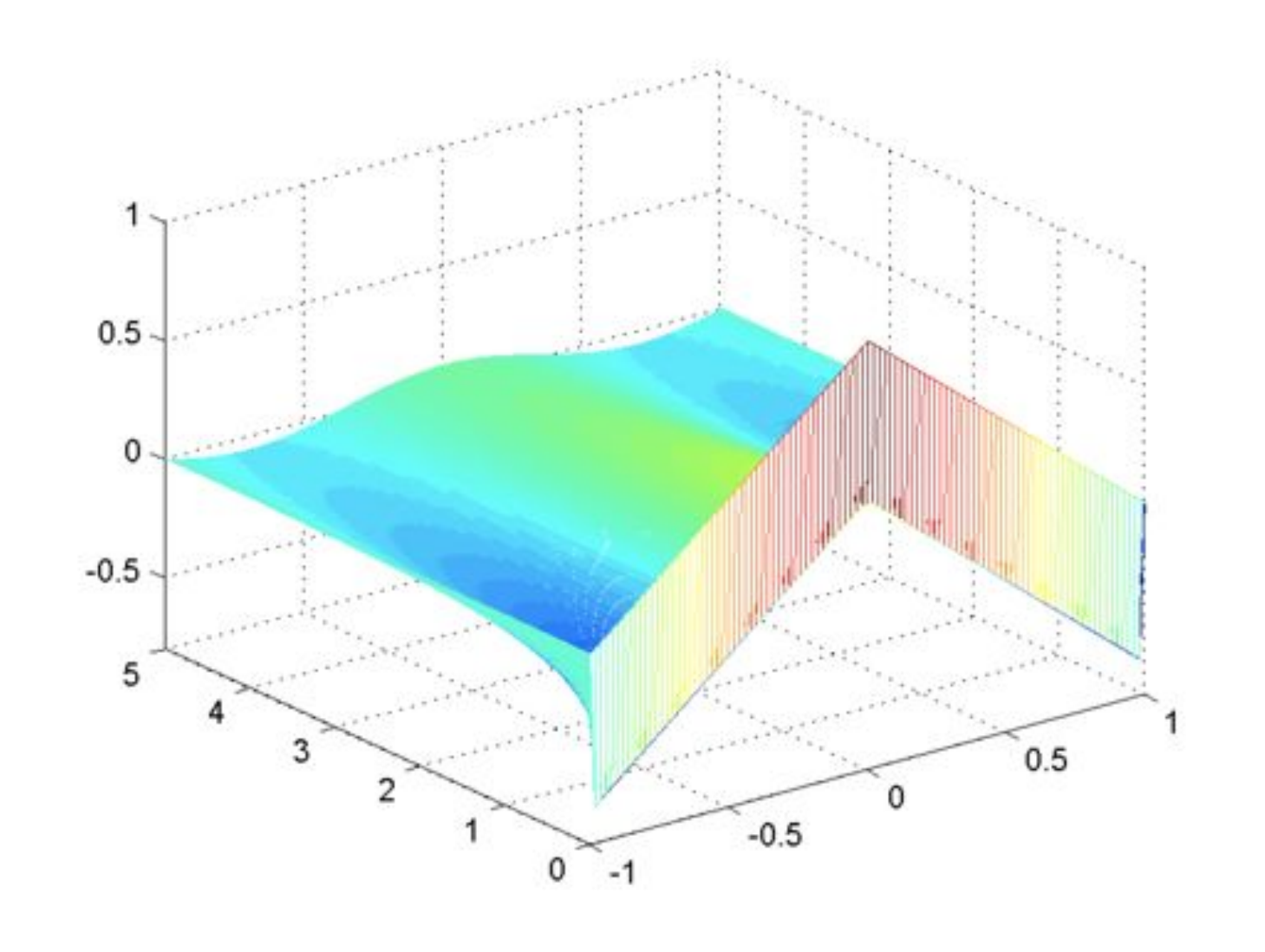}\;
\includegraphics[scale=0.4]{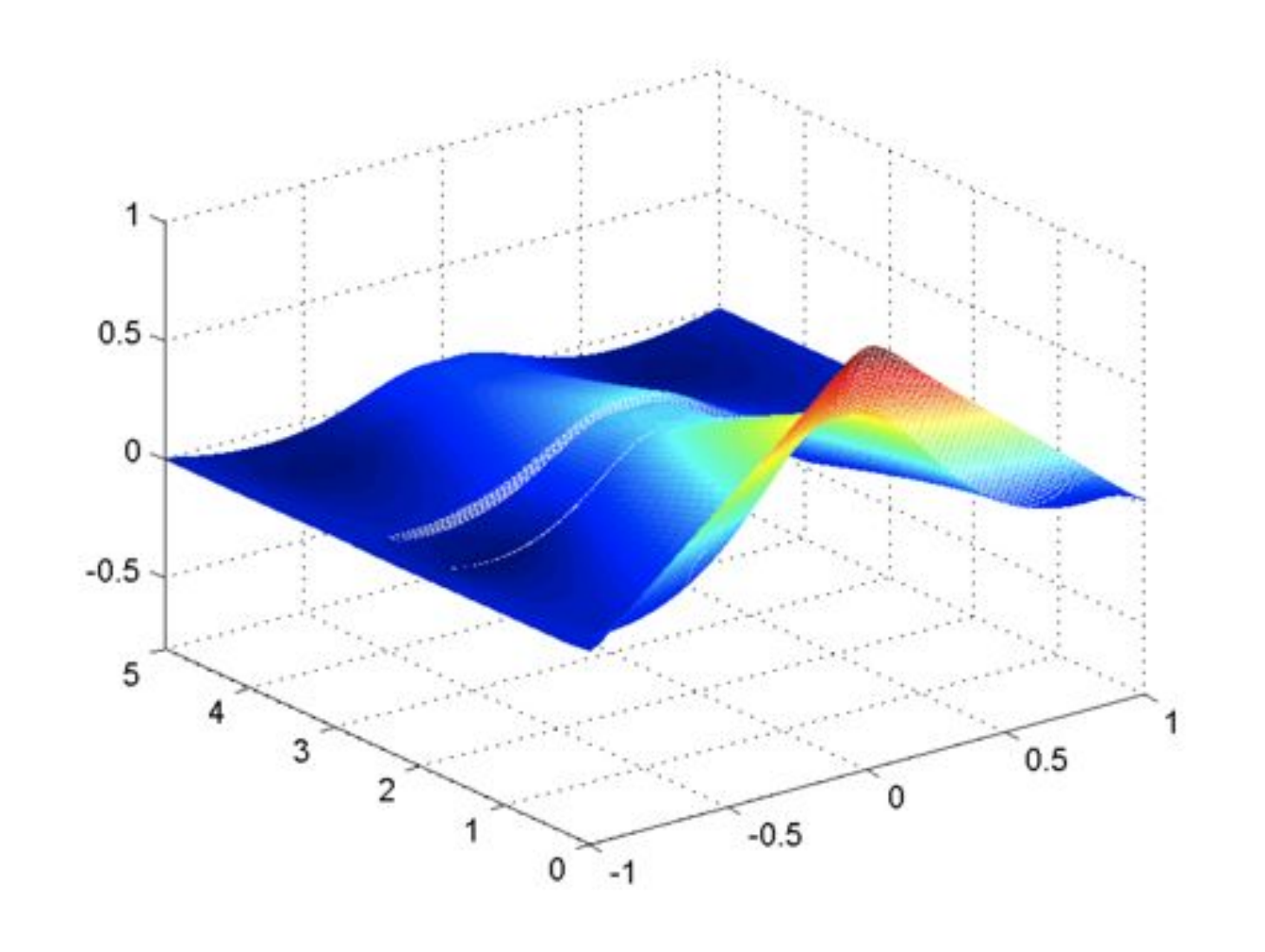}
\end{minipage}

\vspace{0.2cm}
\hspace{0.2cm}(c)\hspace{5.5cm}(d)
\caption{Test 2: (a) exact solution  for $\widehat u=0$; (b) Exact solution for $\widehat u=0$ POD (3 basis functions); (c) Approximate optimal solution for LQR-MATLAB; (d) Approximate solution POD (3 basis functions)+ HJB.}
\label{testass}
\end{figure}

\medskip
\paragraph{Test 2: Heat equation with no-smooth intial data}
In this section we change the initial condition with a function which is only Lipschitz continuos: $y_0(x)=1-|x|.$ According to Test 1, we consider the same parameters. (see Figure \ref{testass}).

\begin{figure}[htbp]
\includegraphics[scale=0.25]{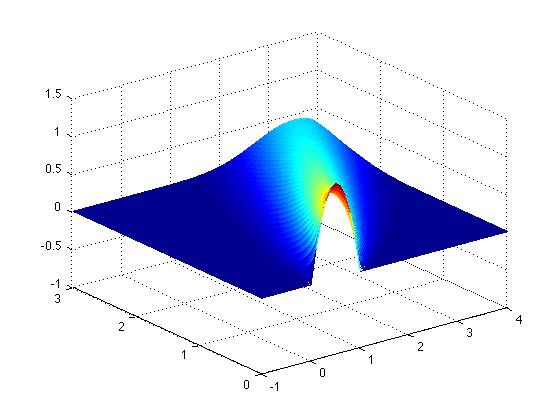}
\hfil\includegraphics[scale=0.25]{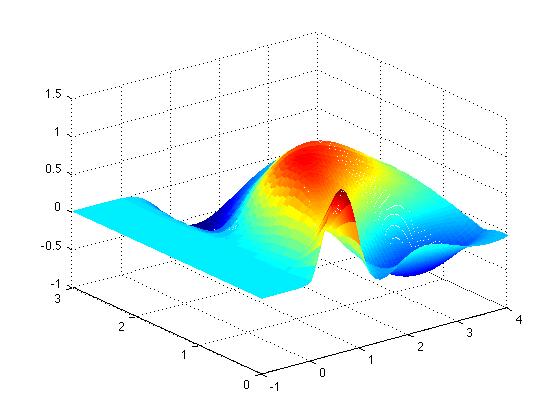}
\caption{Test 3: Solution $\widehat{y}$ on the left, approximate solution on the right with POD (4 basis functions)}
\label{dt_noad}
\end{figure}

\noindent
Riccati's equation has been solved by a MATLAB LQR routine. Thus, we have used  the solution given by this routine as the correct solution  in order to compare the errors in $L^1$ and $L^2$ norm between the reduced Riccati's equation and our approach based on the reduced HJB equation. Since we do not have any information, the snapshots are computed for $\widehat u=0.$ This is only a guess,  but in the parabolic case fits well due to the diffusion term.\\
\noindent
\begin{table}[htbp]
\begin{center}
\begin{tabular}{|c|c|c|}
\cline{1-3}
\hline
  & $L^1$ & $L^2$\\
\hline
& &\\
$y^{LQR}-y^{POD+LQR}$ & 0.0221 & 0.0172\\
\hline
& &\\
$y^{LQR}-y^{POD+HJB}$ & 0.0204 & 0.0171\\
\hline
\end{tabular}
\end{center}
\vspace{0.2cm}
\caption{Test 2: $L^1$ and $L^2$ errors at time $T$ for the optimal approximate solution.}
\label{table1}
\end{table}
\noindent
As in Test 1, the choice of the snapshots does not effect strongly the approximation due to the asymptotic behavior of the solution. 
The presence of a Lipschitz continuous initial condition has almost no influence on the global error (see Table 1).
\medskip
\paragraph{Test 3: Advection-Diffusion equation}
The advection-diffusion equation needs a different method. We can not use the same $\widehat y$ we had in the parabolic case, mainly  because in Riccati's equation the control is free and is not bounded, on the contrary when we solve an HJB we have to discretize the space of controls. We modified the problem in order to deal with  bang-bang controls.  We get $\widehat y$ in (\ref{lqr_oss}) just plugging in the control $\widehat u\equiv 0$.  We have considered the control space  corresponding only  to three values in $[-1,1]$, then $U=\{-1,0,1\}.$ We first have tried to get a controlled solution, without any adaptive method and, as expected, we obtained  a bad approximation (see Figure \ref{dt_noad}). 
From Figure \ref{dt_noad} it's clear that POD with four basis functions is not able to catch the behavior of the dynamics, so we have applied our adaptive method.\\
\begin{figure}
\includegraphics[scale=0.25]{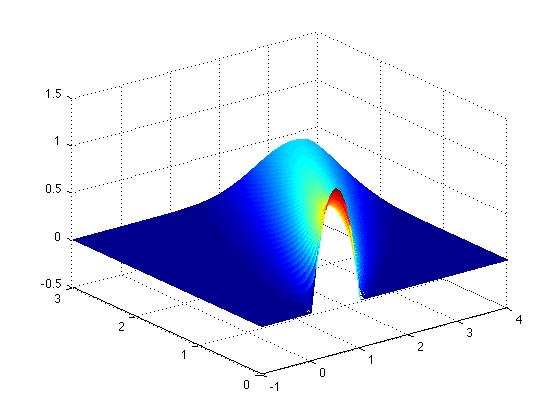}\hfil\includegraphics[scale=0.25]{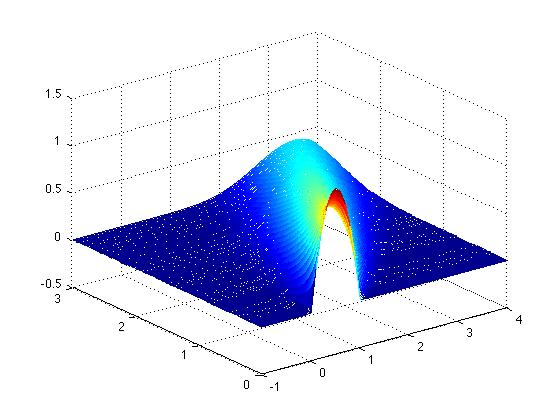}
\caption{Test 3: Solution for $\widehat u\equiv 0$ (left), approximate optimal solution (right).}
\label{dt_ad}
\end{figure}
\noindent
We have consider: $T=3, \Delta x=0.1, \Delta t=0.008$, $a=-1$, $b=4$, $R=0.01.$ According to our algorithm, the time interval $[0,3]$ was divided into $[0,0.744]\cup[0.744, 1.496]\cup[1.496,3].$ As we can see our last interval is bigger than the others, this is due to the diffusion term (see Figure \ref{dt_ad}). The $L^2-$error is 0.0761, and the computation of the optimal solution via HJB  has required about  six hours of CPU time.
In Figure 4 we compare the exact solution with the numerical solution based on a POD representation. Note that, in this case, the choice of only 4 basis functions  for the whole interval $[0,T]$ gives a very poor result due to the presence of the advection term. Looking at Figure 5 one can see the improvement of our adaptive technique which takes always 4 basis functions in each sub-interval.\\
\noindent
In order to check the quality of our approximation we have computed the numerical residual, defined as: 
$$\mathcal{R}(y)=\|y_s(x,s)-\varepsilon y_{xx}(x,s)+cy_x(x,s)-u(s)\|.$$
The residual for the solution of the control problem computed without our adaptive technique is 1.1, whereas the residual for the adaptive method is $2*10^{-2}$. As expected from the pictures, there is a big difference between these two value.

\begin{figure}
\includegraphics[scale=0.25]{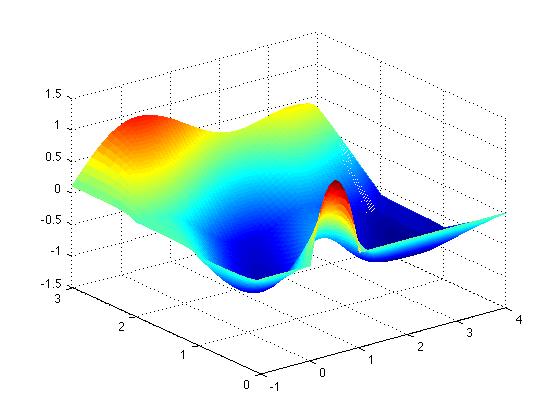}\hfil\includegraphics[scale=0.25]{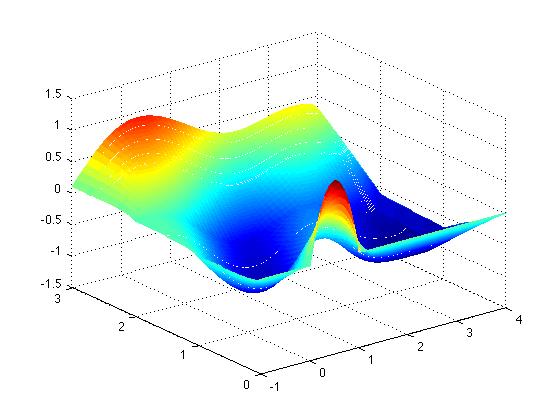}
\caption{Test 4: Solution for $\widehat u$ (left), approximate optimal solution (right).}
\label{dt_ad4}
\end{figure}

\paragraph{Test 4: Advection-Diffusion equation}
In this test we take a different $\widehat y$, namely  the solution of (\ref{eq:test_cont}) corresponding to the control
$$\widehat u(t)=\left\{\begin{array}{l}
-1\quad 0\leq t< 1\\
0\quad 1\leq t< 2\\
1\quad 2\leq t\leq 3.
\end{array}
\right.$$
We want to emphasize we can obtain nice results when the space of controls has few element. The parameters were the same used in Test 3. The $L^2-$error is 0.09, and the time was the same we had in Test 3. In Figure \ref{dt_ad4} we can see our approximation.
In Figure 6 one can see that the adaptive technique can also deal with discontinuous controls.\\
\noindent

In this test, the residual for the solution of the control problem without our adaptive technique is 2, whereas the residual for the adaptive method is $3*10^{-2}$. Again, the residual shows the higher accuracy of the adaptive routine.

%

\section{Conclusions}

As we have discussed,  a reasonable coupling between POD and HJB equation 
can produce feedback controls for infinite dimensional problem. For advection dominated 
equations that simple idea has to be implemented in a clever way to be successful.
It particular, the application of an adaptive technique is crucial to obtain accurate approximations with a low  number of POD basis functions. This is still an essential requirement when dealing with the Dynamic Programming approach, which suffers from the curse-of-dimensionality although recent developments in the methods used for HJB equations will allow to increase this bound in the next future (for example by applying patchy techniques).

Another important point is the discretization of the control space. In our examples, the number of optimal control is rather limited and this will be enough for problems which have a bang-bang structure for optimal controls. In general, we will need also an approximation of the control space via reduced basis methods. This point as well as a more detailed analysis of the procedure outlined in this paper will be addressed in our future work.





\end{document}